\documentclass[12pt]{amsart}
\usepackage[centertags]{amsmath}
\usepackage{amscd,latexsym,amssymb} 
\usepackage{euscript,mathrsfs} 
\newcommand{\B}[1]{{\mathbb #1}}

\newcommand{\Z}{\B Z} 
\newcommand{\R}{\B R}  
\newcommand{\T}{\B T}

\newtheorem{theorem}[subsection]{Theorem}

\newtheorem{cory}[subsection]{Corollary} 
\newtheorem{lemma}[subsection]{Lemma} 
\newtheorem{proposition}[subsection]{Proposition} 
\newtheorem{prop}[subsection]{Proposition}

\theoremstyle{definition} 
\newtheorem{definition}[subsection]{Definition} 
\newtheorem{example}[subsection]{Example} 
\newtheorem{df}[subsection]{Definition} 
\newtheorem{constructions}[subsection]{Constructions} 
\theoremstyle{remark} 
 
\newtheorem{question}[subsection]{Question}

\numberwithin{figure}{section} 
\numberwithin{table}{section} 
\numberwithin{equation}{section}

\def\m{\medskip}

\newcommand{\Om}{{\Omega}} 
\newcommand{\om}{{\omega}}

\newcommand{\Ga}{{\Gamma}}

\newcommand{\la}{{\lambda}}

\newcommand{\Si}{{\Sigma}} 
 
\newcommand{\gf}{\varphi}

\def\la{\langle} 
\def\ra{\rangle}
\def\pg{p_{\Ga}}
\def\wh{\widehat}
\def\ov{\overline}

\newcommand{\map}[1]{\stackrel {#1}\longrightarrow}

\newcommand{\Mo}{(M,\omega )} 
 
\newcommand\pt{\operatorname{pt}}

\newcommand\ab{\operatorname{ab}} 
\newcommand\Diff{\operatorname{Diff}} 

\newcommand\Hom{\operatorname{Hom}}

\newcommand\Ker{\operatorname{Ker}}

\newcommand\supp{\operatorname{supp}} 
\newcommand\rank{\operatorname{rank}}

\newcommand\theoref{Theorem~\ref}
\newcommand\lemref{Lemma~\ref}
\newcommand\propref{Proposition~\ref}
\newcommand\corref{Corollary~\ref}

\newcommand\coryref{\corref}

\newcommand\AAA{\mathcal A} 
\newcommand\BBB{\mathcal B}

\long\def\forget#1\forgotten{} %
\begin{document}

\title[Symplectically aspherical abelian groups]
{On fundamental groups of 
symplectically aspherical manifolds II: abelian groups}

\author{J. K\c EDRA}
\author{Yu. RUDYAK}
\author{A. TRALLE}

\address{J. K\c edra, Mathematical Sciences, University of Aberdeen,
Meston Building, King's College, Aberdeen AB24 3UE, Scotland, UK\hfill\hfill\newline
{and}\hfill\hfill\newline
Institute of Mathematics, University of Szczecin,
Wielkopolska 15, 70451 Szczecin, Poland}
\email{kedra@maths.abdn.ac.uk}
\urladdr{http://www.maths.abdn.ac.uk/\~{}kedra}

\address {Yu. Rudyak, Department of Mathematics, University of Florida,
358 Little Hall, Gainesville, FL 32601, USA}
\email{rudyak@math.ufl.edu} 
\urladdr{http://www.math.ufl.edu/\~{}rudyak}

\address{A. Tralle, Department of Mathematics, University of Warmia and 
Mazury,
10561 Olsztyn, Poland}
\email{tralle@matman.uwm.edu.pl}  
\urladdr{http://wmii.uwm.edu.pl/\~{}tralle} 





\begin{abstract}  
We describe all abelian groups which can appear as the 
fundamental groups of closed symplectically aspherical manifolds. 
The proofs use the theory of symplectic Lefschetz fibrations.
\end{abstract} 

\maketitle

\section{Introduction}\label{S:intro} 
 
The present paper is a continuation of \cite{IKRT}. 
Recall that a symplectic form $\omega$ on a smooth  manifold $M$ is called
{\em symplectically aspherical} if 
$$
\int_{S^2}f^*\omega=0
$$
for any continuous map $f: S^2\rightarrow M$. 
This is expressed equivalently by saying that the cohomology class 
$[\omega]\in H^2(M;\mathbb{R})$ vanishes on the image of the Hurewicz map 
$h: \pi_2(M)\rightarrow H_2(M;\mathbb{R})$. A {\em symplectically aspherical manifold} is defined as a connected manifold that admits a symplectically aspherical form.

Gompf asked a question about the topology of such manifolds 
in \cite{G2}. The importance of this class of symplectic manifolds
comes from the Floer theory  which is much simpler in the
symplectically aspherical case \cite{F, S}. 
In \cite{IKRT} and in this article we are interested in the 
following question (cf. \cite{G2}).

\begin{question}\label{q:real}
What groups can be realized as fundamental groups of closed
symplectically aspherical manifolds?
\end{question}

In the sequel we call these groups {\it symplectically aspherical}.
One of our main results is the classification of finitely generated
abelian symplectically aspherical groups.

\begin{theorem}\label{T:abel}
A finitely generated abelian group $\Gamma $ is symplectically aspherical
if and only if either $\Gamma \cong \mathbb Z^2$ or $\rank(\Ga)\geq 4$.
\end{theorem}

\m
Realization Question~\ref{q:real} can be thought of in the general
context of possible restrictions 
which  an additional geometric structure on a given manifold imposes on its
algebraic topology and, in particular, on its fundamental group. 
For example, we could mention the Hopf question on fundamental groups of
3-manifolds or Serre's question about fundamental groups of 
projective varieties (see Problem on page 10 in \cite{Se} or
\cite[Open Problem 1.17]{ABCKT}). 
The same question for compact K\"ahler manifolds is still 
unsolved (see \cite{ABCKT} for an exposition). 

\m In \cite{G2} Gompf showed that any finitely presentable group can
be realized as the fundamental group of a closed symplectic manifold.
On the other hand, this general realization result fails under the condition of 
symplectic asphericity: for instance, finite groups and $\Z$ are not  
symplectically aspherical, cf. \cite{IKRT}. 

\begin{example}\label{E:3dim} 
Let $\Ga$ be a (finitely generated) group of real cohomological
dimension 3, i.e. $H^3(\Ga;\R)\ne 0$  while $H^i(\Ga;\R)=0$ for $i>3$. 
Then $\Ga$ is not  symplectically aspherical. 
Indeed, if $\Mo$ is a closed symplectically aspherical manifold and 
$\pi_1(M)=\Ga$ then the class $[\om]$ lies 
in the image of the homomorphism $c^*:H^2(K(\Ga,1);\R)\to H^2(M;\R)$, 
induced by the classifying map $c: M \to K(\Ga,1)$. In other words, $[\omega]=c^*(a)$ for some 
$a\in H\sp 2(\Ga;\R)$.
Since $\Ga$ is three dimensional, we conclude that $[\om]^2=c^*(a^2) =0$, and hence
$\dim M = 2$. Therefore $M$ is an aspherical closed oriented  
surface, and thus $H^3(\Ga;\R)=H\sp 3(M;\R)=0$. 
\end{example} 

\m We see that, in particular, no finitely generated abelian group of rank 3 are  
symplectically aspherical. Hence our Theorem \ref{T:abel} states that the obvious necessary 
condition is also sufficient for a finitely generated abelian
group to be symplectically aspherical.

\m
Moreover, we answer the question from \cite{IKRT}, motivated by
Gompf \cite{G2}, about the relation between two classes of symplectically 
aspherical groups. 
Recall that, the class $\AAA$ consists of groups $\Gamma$ realizable as
$\pi_1(M)$, where $M$  s symplectically aspherical with $\pi_2(M)=0$, while 
the class $\BBB$ consists of symplectically aspherical groups realizable
as $\pi_1(M)$ with  $\pi_2(M)\not=0$. 
It is easy to see that the group $\Z^2$ belongs to $\AAA$ and does not
belong to $\BBB$,  while in \cite{IKRT} we asked whether $\BBB
\subset \AAA$.  
In this paper we show that $\BBB\not\subset\AAA$, namely, $\Z\sp 4 \oplus \Z/2\in \BBB \setminus \AAA$(see \propref{prop:pi2}). This  phenomenon deserves further investigation,
since this gives a  non-realizability result, which may reveal some
new essentially symplectic (non-topological) properties.

\m
In this work we focus on 4-dimensional symplectic Lefschetz fibrations
and provide conditions implying the symplectic asphericity of the
total space. The use of this technique proved to be very effective in
many other contexts. In fact, any
closed 4-dimensional  symplectic manifold admits a Lefschetz pencil, and, 
therefore, becomes a symplectic Lefschetz fibration after blow-up at
several points \cite{D1, GS}. 
Using the Donaldson hyperplane section theorem [D2], one can reduce
the realizability  problem to the 4-dimensional case 
(see \cite[Proposition 2.2]{IKRT}
for details). Explicit construction of Lefschetz 
fibrations with given fundamental groups was given by
Amoros, Bogomolov, Katzarkov and Pantev \cite{ABKP}. 
Our results can be regarded as solutions of a similar construction 
problem under additional restriction of symplectic asphericity.    

\m
Throughout the paper $\Sigma_g$ denotes the closed orientable surface of genus $g$ and 
$\pi_g$ denotes the fundamental group of $\Sigma_g$.

\m {\bf Acknowledgments:} The second author was supported by NSF,
grant 0406311. The third author was supported by the Ministry of Education
and Science of Poland, research grant
no. 1P03A 03330.

We are also grateful to Ivan Smith for answering our questions.


\section{Symplectic Lefschetz fibrations}\label{s:slf}


\begin{df}[\cite{ABKP, D1, GS}] Let $X$ be a compact, connected, oriented, 
smooth 4-manifold, possibly with boundary. A {\em Lefschetz fibration structure} on $X$ is 
a surjective map $f: X \to \Si$ where $\Si$ is a compact, connected, oriented surface and 
$f^{-1}(\partial \Si)=\partial X$. Furthermore, the following is required:

\begin{itemize}
\item the set $\{q_1, \ldots, q_n\}$ of critical 
points of $f$ is finite;
\item $f(q_i)\neq f(q_j)$ for $i\neq j$; 
\item if $b\in \Si_{g}$ is a regular value of $f$ then $f^{-1}(b)$ is a closed 
connected orientable surface;
\item let $F_i=f^{-1}(f(q_i))$.  Then there exists a complex chart $\gf: U_i\to  
\B C^2$ with $q_i\in U_i\subset X$ and a complex chart $\psi_i:V_i\to \B C$ with 
$f(q_i)\in V_i \subset f(U)\subset  \Si$ such that $\psi_if\gf_i^{-1}: \B C^2 \to 
\B C$ has the form  $(x,y) \mapsto x^2+y^2$.
Here we require that both complex charts preserve the orientations.
\end{itemize}
\end{df}

Clearly, the manifolds $f^{-1}(b_1)$ and $f^{-1}(b_2$) are diffeomorphic for any 
two regular values $b_1, b_2$ of $f$. We define the {\it regular fiber} of $f$ to 
be the (diffeomorphism equivalence class of the) manifold $f^{-1}(b)$ where $b$ 
is a regular value of $f$. Non-regular fibers are called {\it singular fibers}.

\m
When we write ``a Lefschetz fibration $F\to X \to B$'', it means that 
$F$ is the regular fiber.

\begin{constructions}\label{constr}

(a) Let $D$ be a closed disc $D\subset \B C$ about the origin and $F\to X \stackrel{f}\to D$ be a Lefschetz fibration over $D$ with one singular point over $0\in D$. The restriction of $f$ onto $\partial D$ is the locally trivial bundle that gives us and is  characterized by a self-diffeomorphism  $T: F \to F$ that is either the identity map or (isotopic to) the Dehn twist along a certain simple closed curve $C$ on $F$.  In other words, the monodromy about the origin of the disc is the (isotopy class of the) Dehn twist $T: F \to F$ along $C$. This curve $C$ is a {\em vanishing cycle} corresponding to $T$, and the singular fiber of $f$ is homeomorphic to $F/C$.

\m 
(b) In order to construct a Lefschetz fibration
$F\to X\stackrel{f}\to S^2$ over the sphere
one proceeds as follows. Following item (a), take $n$ Lefschetz
fibrations $X_1, \dots ,X_n$ over discs $D_1, \ldots, D_n$ with
monodromies $T_1,\dots ,T_n$, respectively. 
Take the boundary connected sum $D'$ of these discs and consider the corresponding 
fiber sum of the Lefschetz fibrations $X_i$. In this way we  
obtain a Lefschetz fibration
$X'\to D'$ over the disc $D'$ with monodromy over $\partial D'$ equal to
the product $T_1\cdots T_n$.
Assuming that this product is isotopic to the identity and
choosing an isotopy, one extends the above fibration $X'\to D'$ 
to a fibration $X_0\to D_0$ over a larger disc $D_0\supset D'$ so that the new one 
has trivial monodromy over the boundary. Finally, 
this fibration over $D_0$ can be extended trivially to a 
Lefschetz fibration over the sphere.
\end{constructions}

\m Each Dehn twist $T_i$ in \ref{constr}(b) can be chosen so that it is supported in a small 
tubular neighborhood of its vanishing cycle. Let 
$p\in F \setminus (\supp T_1 \cup\dots\cup \supp T_n)$ be
a basepoint.

\begin{lemma}\label{L:section}
Suppose that the product $T_1\cdots T_n$ 
represents the neutral element in the mapping
class group $\pi_0(\Diff(F,p))$ {\rm(}in other words, the isotopy
from $T_1\cdots  T_n$ to the identity map preserves
$p\in F${\rm )}. Then the Lefschetz fibration constructed in {\rm \ref{constr}(b) }
admits a section. Furthermore, $\pi_1(X)$ is isomorphic 
to the quotient group $\pi_1(F)/N$ where $N$ is the normal subgroup of $\pi_1(F)$ generated by 
vanishing cycles $C_1,\ldots, C_n$ corresponding to Dehn twists $T_1, \ldots, T_n$, respectively. 
\end{lemma}

\begin{proof}
First observe that each fibration 
$F\to X_i \to D_i$ over the small disc
about a critical value contains the trivial subbundle
$(F \setminus \supp T_i) \times D_i$. Let 
$s_i:D_i\to (F \setminus \supp T_i) \times D_i\subset X_i$ 
be the section given by $s_i(x):= (p,x)$.
Together these sections yield a section
$s':D'\to (F\setminus(\supp T_1 \cup \dots \cup \supp T_n ))\times D'\subset X'$,
$s'(x):=(p,x)$.
Since the isotopy from $T_1\cdots T_n$ to the identity map preserves
the basepoint $p$, we conclude that the section $s'$ extends to a section
$S^2\to X$.

The last claim is proved in \cite[Lemma 2.3]{ABKP}.
\end{proof}

\m A {\it symplectic Lefschetz fibration} is defined to be a Lefschetz fibration 
with a symplectic form $\omega$ on the total space whose restriction on each of the  
fibers is non-degenerate, see \cite{ABKP, GS} for details.

\m
A proof of the following \theoref{theor:amoros} (without item (iv)) is contained in 
\cite[Theorem A]{ABKP}. A weak version of this theorem can also be deduced by combining the theorem 
of Donaldson \cite{D1, ADK} that every symplectic 4-manifold admits a structure 
of a Lefschetz fibration, and the theorem of Gompf \cite{G1} on realizability of 
any finitely presentable group as the fundamental group of closed symplectic 
4-manifold.

\m Recall that a group $G$ is called {\em finitely presentable} if there exists an 
epimorphism $f:F \to G$ where $F$ is a free group of 
finitely many free generators and $\Ker f$ is a normal closure of a finitely
generated subgroup of $F$. Following \cite{ABKP}, we define a 
finite presentation of a group $G$ to be an epimorphism $a: A \to G$ where $A$ 
is a finitely presentable group and $\Ker a$ is a normal closure of a finitely
generated subgroup of $A$.

\begin{theorem}\label{theor:amoros}
Let $p_{\Ga}:\pi_g\to \Ga$ be a finite presentation of a group $\Ga$. Then there 
exists an 
epimorphism $p_{h,g}: \pi_h\to \pi_g$ for some $h\ge g$ and a symplectic 
Lefschetz fibration
$$
\CD
\Si_h @>i>> X @>f>> S^2
\endCD
$$
with the following properties:
\begin{enumerate}
\item[(i)]$\pi_1(X)=\Gamma$;
\item[(ii)] the homomorphism $i_*: \pi_h \to \pi_1(X)$ coincides with $p_{\Ga}p_{h,g}$;
\item[(iii)] the homomorphism $p_{h,g}$ is induced by a map $\Si_h\to \Si_g$ of non-zero degree;
\item[(iv)] the map $($fibration$)$ $f$ has a section.
\end{enumerate}
\end{theorem}

\begin{proof}The first three items are proved in \cite[Theorem A]{ABKP}. Roughly speaking, the idea of the proof looks as follows. The authors construct a map $p_{h,g}:\pi_h \to \pi_g$ such that the map $p_{\Ga}p_{h,g}: \pi_h \to \Ga$ is a finite presentation of $\Ga$. Moreover, the kernel of $p_{\Ga}p_{h,g}$ is generated as a normal subgroup by homotopy classes of simple closed curves $C_1, \ldots, C_n$. Furthermore, the Dehn twists $T_i$ along $C_i$, $i=1, \ldots, n$ satisfy the conditions of \lemref{L:section}. Now, we construct a Lefschetz fibration as in \ref{constr}(b). This Lefschetz fibration turns out to be symplectic in view of \cite[Theorem 10.2.18]{GS} or \cite[Proposition 2.3]{ABKP}.

Because of what we said above and \lemref{L:section}, we get the proof of (i), (ii) and (iv). 

Finally, the map $p_{h,g}$ is induced by a composition $\Si_h \to \Si_e\to \Si_g$, where $\Si_e$ is obtained from $\Si_g$ by adding handles (and the map $\Si_e\to \Si_g$ collapses these handles) and the map $\Si_h \to \Si_e$ is a finite ramified covering. This implies (iii).
\end{proof}

\m
Recall that, given a map $f:X \to Y$ of spaces, a cohomology class $a\in H^k(X)$ 
is said to be {\em totally non-cohomologous to zero}, abbreviated as  TNCZ, if 
$i_y^*(a)\neq 0$ for each inclusion $i_y:f^{-1}(y)\subset X$ where $y$ runs over 
all points of $Y$.

\begin{theorem}\label{theor:gs}
Let $\omega_{\Si_g}$ be a symplectic form on $\Si_g$ and 
$$
\CD
F @>i>>  X @>f>> \Si_g
\endCD
$$ 
be a symplectic Lefschetz fibration. Let $\Omega\in 
H^2(X;\R)$ be a  TNCZ class. Then there exists $C\in \R$ such 
that the cohomology class $\Omega+Cf^*[\omega_{\Si_g}]$ contains a symplectic 
form on $X$.
\end{theorem}

\begin{proof} This is actually proved in \cite[Theorem 10.2.18]{GS}. 
Namely, without loss of generality we can assume that $\la \Omega, [F]\ra >0$. 
Then $\la \zeta, [F]\ra >0$ for each closed form $\zeta\in \Omega$. Let $\eta$ 
be the form described in the proof of Theorem 10.2.18 in \cite{GS}. 
Then $\eta \in\Omega$, while $\eta +Cf^*\omega_{\Si_g}$ is a symplectic form on 
$X$, see \cite[Proposition 10.2.20]{GS}.
\end{proof}


\section{Symplectically aspherical Lefschetz fibrations}\label{S:salf}


Given a space $X$, we call a non-zero cohomology class $a\in H^k(X;G)$  {\em aspherical}, 
if $\langle a,f_*[S^k]\rangle=0$ for any continuous map $f:S^k\to X$. 

\begin{proposition}\label{P:fac} 
Let $p_{\Ga}: \pi_g\to \Ga $ be a finite presentation of a group $\Ga$ such 
that 
$$
\pg^*: H^2(\Ga;\R)\to H^2(\pi_g;\R)
$$ is a non-zero homomorphism.
Then there exists a symplectic Lefschetz fibration 
$\Sigma_h\stackrel{i}\to X \stackrel{f}\to S^2$ with $\pi_1(X)=\Ga$ 
and TNCZ aspherical class $\Om\in H^2(X)$. Furthermore, this Lefschetz fibration admits a section.
\end{proposition} 

\begin{proof}
Because of \theoref{theor:amoros}, there exists an epimorphism 
$$
p_{h,g}:\pi_h \to \pi_g
$$
and a symplectic Lefschetz fibration $\Si_h\stackrel{i}\to X \stackrel{f}\to 
S^2$ admitting a section and such that $\pi_1(X)=\Ga$ and $i_* = p_{\Ga}\circ p_{h,g}$. Moreover, since $p_{h,g}$ is induced by a map $\Sigma_h \to \Sigma_g$ of non-zero degree, we conclude that  
\begin{equation}\label{eq:phg}
p_{h,g}^*: H^2(\pi_g; \R) \rightarrow H^2(\pi_h;\R)
\end{equation}
is an isomorphism. 

Let $c: X\to K(\Ga,1)$ be a map that induces an isomorphism of fundamental 
groups.  Let $i: \Si_h \to X$ be the inclusion of a regular fiber. We have the 
commutative diagram
$$
\CD
\pi_h @=\pi_1(\Si_h)@>i_*>>\pi_1(X) @>c_*>\cong >\pi_1(K(\Ga,1))\\
@| @VVp_{h,g}V @| @| \\
\pi_h @>p_{h,g}>> \pi_g @>p_{\Ga}>> \Ga @=\Ga
\endCD
$$
Since the homomorphism \eqref{eq:phg} is an isomorphism, we conclude that the homomorphism $i^*\circ c^*:H^2(\Ga; \R) \to H^2(\pi_1(\Si_h);\R)$ is nontrivial because so is $p_{\Ga}^*: H^2(\Ga; \R)\to H^2(\pi_g; \R)$. Take any $a\in H^2(\Ga;\R)$ with $c^*(a)\ne 0$ and put $\Omega=c^*(a)$. Then $\langle \Omega, i_*[\Si_h])\rangle \ne 0$ because $i^*( c^*(a))\ne 0$.

Finally, if $j: \Si \to X$ is the inclusion of a singular fiber then $\langle \Omega, j_*(\Si)\rangle \ne 0$ because singular fibers are homologous in $X$ to regular fibers.
\end{proof}

\begin{lemma}\label{L:fac} Let $F\stackrel{i}\to X  \stackrel{f}\to \Si_g$ be a symplectic 
Lefschetz fibration over the base of genus $g>0$. Suppose that  
$X$ admits a TNCZ aspherical class $\Om\in H^2(X;\R)$. Then 
$X$ is symplectically aspherical.  
\end{lemma} 

\begin{proof} 
Because of \theoref{theor:gs}, there exists a symplectic 
structure $\om$ on $X$ whose cohomology class is of the 
form $\Om + C f^*[\om_{\Si_g}]$, for some constant $C\in \B R$.  
Clearly, it is an aspherical class. 
\end{proof} 

\begin{prop}\label{p:sum}Let $F\to X \to \Si_g$ be a symplectic Lefschetz fibration 
that admits a TNCZ aspherical class $\Omega\in H^2(X;\R)$. Let $Y =F\times \Si_h$ with the product
symplectic structure. If $g+h>0$ then the Gompf symplectic
fiber sum $X \#_F Y$ is symplectically aspherical.
\end{prop}

\begin{proof}
It is clear that $F \to X \#_F Y\to \Sigma_g \# \Sigma_h = \Sigma_{g+h}$
is a symplectic Lefschetz fibration.
Furthermore, the retraction $\Si_h \to D^2$ yields the retraction 
$\Si_h\times F \to D^2\times F$.
Now, the degree one map
$\Sigma_g\# \Sigma_h \to \Sigma_g$ (collapsing the summand $\Sigma_h \setminus D^2$
to the point) yields the following commutative diagram:
$$
\CD
X\#_{F} Y @>\gf>> X\\
@VVV @VVV\\
\Si_g\#\Si_h @>>> \Si_g.
\endCD
$$
Clearly, the class $\gf^*\Omega$ is aspherical. Furthermore,  $\gf^*\Omega$ is TNCZ 
since $\gf$ is a fiber map, and so $\la \gf^*\Omega, [F]\ra=\la \Omega, [F]\ra$. Now 
the statement follows from \lemref{L:fac}. 
\end{proof}


\section{Fundamental groups of fiber sums}\label{S:fibersum}


\begin{definition}
A {\em short surjectivity diagram} is a diagram of groups and homomorphism of 
the form
\begin{equation}\label{eq:ssd}
\CD
1 @>>> A @ >j>> B @>>> G @>>> 1\\
@. @V\gf VV @VV \psi V @| @.\\
1 @>>> P @ >>> Q @>>> G @>>> 1
\endCD
\end{equation}
where $\gf$ is an epimorphism and the rows are exact.
\end{definition}

\begin{prop}\label{p:ssd}
In the short surjectivity diagram \eqref{eq:ssd} the map $\psi$ yields an 
isomorphism $B/j(\Ker \gf)\to Q$. Furthermore, if the top exact sequence splits 
then so does the bottom sequence.
\end{prop}

\begin{proof} The first claim is clear by diagram chasing, cf. \cite[Lemma II.3.2]{M}. The 
second claim follows from the first one.
\end{proof}

\begin{prop}\label{p:circle}
Let $F\map{i} X \map{f} S^2$ be a symplectic Lefschetz fibration that admits a section. Then, for 
every regular fiber $F$, the inclusion $X\setminus F \to X$ induces an isomorphism of 
fundamental groups.  
\end{prop}

\begin{proof}
Let $D=D^2\subset S^2$ be a small closed disk centered at $f(F)$, and let 
$U=f^{-1}(D)$. Then $\partial U= S^1\times \pi_1(F)$. Take the section $s: S^2 \to X$ and restrict 
it onto $\partial (S^2\setminus D)=\partial D$. Now, we regard the (pointed) 
homotopy class $\alpha$ of $s: \partial D \to \partial U$ as an element of $\pi_1(\partial U)$, 
and it is clear that $\pi_1(\partial U)\cong \Z \times \pi_1(F)$ where $\Z$ is a subgroup 
generated by $\alpha$. 

Because of the Seifert--van Kampen theorem, we have an isomorphism
\[
\pi_1(X) \cong\pi_1(X\setminus F)*_{\pi_1(\partial U)}\pi_1(U)= \pi_1(X\setminus 
U)*_{\Z \times \pi_1(U)}\pi_1(U).
\]
Since the fibration has the global section, we conclude that the subgroup $\Z=\{\alpha\}$ maps to zero
under both inclusions $\partial U \subset X \setminus U$ and $\partial U \subset\overline U$. Thus
\[
\pi_1(X)=\pi_1(X\setminus F)*_{\pi_1(U)}\pi_1(U)=\pi_1(X\setminus F)*_{\pi_1(F)}\pi_1(F)=\pi_1(X\setminus F).
\]
The proof is completed.
\end{proof}

\begin{lemma}\label{L:pi1} 
Let $F \stackrel{i}\to X \to S^2$ be a symplectic Lefschetz fibration that admits a section. 
Let 
\begin{equation}\label{bundle}
F\to Y\stackrel{p}\to \Si_e
\end{equation}
be a surface bundle with $e>1$ that  admits a section. 
Then there exists a short surjectivity diagram
\begin{equation*} 
\CD
1 @>>> \pi_1(F) @>>> \pi_1(Y) @>>> \pi_e @>>>  1\\ 
@. @Vi_*VV @VVV @| @. \\
1 @>>> \pi_1(X) @>>>\pi_1(X\#_F Y) @>>> \pi_e  @>>>1
\endCD
\end{equation*}
where the upper row is the segment of the homotopy exact sequence of the bundle $\eqref{bundle}$.
\end{lemma} 

\begin{proof} 
Let $f$ be the genus of $F$ and 
$$
\left\langle a_1,b_1,\dots,a_f,b_f\,\, \bigg\vert 
\prod\,[a_i,b_i]\right\rangle
$$
be the standard presentation of $\pi_f$. Notice that $i_*$ is surjective because of the exactness of the sequence 
\[
\pi_1(F) \to \pi_1(E) \to \pi_1(B)
\]
for every Lefschetz fibration $F \to E\to B$, see \cite[Proposition 8.1.9 and p. 510]{GS}. 
So, $\pi_1(X)$ has a presentation
$$ 
\left\langle \wh a_1,\wh b_1,\dots,\wh a_f,\wh b_f \bigg\vert \prod[\wh a_i,\wh b_i],\wh R\right\rangle
$$ 
where $\wh R$ is a finite set of words and $\wh a_k=i_*a_k,\wh b_k=i_*b_k, k=1, \ldots, f$.

Let $x_1,y_1,\dots x_e,y_e$ be the standard system of meridians and parallels on $\sigma(\Si_e)$, where $\sigma: \Si_e \to Y$ is a section of the bundle \eqref{bundle}. Regarding $x_i$ and $y_i$ as elements of $\pi_1(Y \setminus F)$, we obtain a presentation
$$
\left\langle x_1,y_1,\dots x_e,y_e,a'_1,b'_1,\dots,a'_f,b'_f\,\, \bigg\vert 
\prod\,[a'_i,b'_i],q'_1,\dots,q'_n\right\rangle
$$ 
of $\pi_1(Y\setminus F)$, where $a'_k=j_*a_k,b'_k=j_*b_k, k=1,\ldots, f$. In fact, $\pi_1(Y\setminus F)$ is the semidirect product $\pi_1(\Si_e\setminus \pt)\ltimes \pi_1(F)$.

Take a small disk $D'\subset \Si_e$ about a regular value of $p$ and consider $u'=p^{-1}(D')\cong D'\times F$. The map $\sigma: \partial D' \to \partial U'=\partial D'\times F$ gives us an element $\beta\in \pi_1(\partial U')$. Notice that the image of $\beta$ in $\pi_1(Y\setminus F)$ is equal to $\displaystyle \prod_{i=1}^e[x_i,y_i]$.

In order to perform the fiber sum $X \#_F Y$ we identify (fiberwisely) the neighborhood $U'$ with the neighborhood $U$ from \propref{p:circle}. It turns out to be that, under this identification, $\beta$ coincides with $\alpha$ described in \propref{p:circle}. This is true because both sections $s: \partial D \to X$ and $\sigma: \partial D'\to Y$ extend to the whole discs $D$ and $D'$, respectively. 

Consider the group $\Z=\{\alpha\}=\{\beta\}$.
Because of the isomorphism $\pi_1(X \setminus F) \cong \pi_1(X) $ from \propref{p:circle} and by the Seifert -- van Kampen theorem, the fundamental group $\pi_1(X \#_F Y)$
of the fiber sum is the amalgamated product 
$$ 
\CD 
\B Z\times \pi_f @>\ov i>> \pi_1(Y\setminus F)\\
@V\ov jVV @VVV\\
\pi_1(X) @>>> \pi_1(X \#_F Y)
\endCD
$$ 
where
\begin{eqnarray*} 
\ov i(\alpha,e_1)= \Pi_{i=1}^{e}[x_i,y_i],  
\quad i(0,a_i) = a'_i, \quad \ov i(0,b_i)=b'_i\\  
\ov j(1,e_1)= e_2, \quad \ov j(0, a_i) = \wh a_i,  
\quad \ov j(0,b_i)=\wh b_i\\  
\end{eqnarray*}
and $e_1, e_2$ are the neutral elements of $\pi_f$ and $\pi_1(X)$, respectively. 
Thus the presentation of $\pi_1(X\#_F Y)$ has the following 
form 
$$ 
\left\la a_1,b_1,\dots,a_f,b_f,x_1,y_1,\dots,x_e,y_e \,\bigg\vert
\Pi[a_i,b_i],\Pi[x_i,y_i], 
q_1,\dots, q_n, R 
\right\ra 
$$ 
where $q_i$ are monomials $q'_i$ with $a'_i$ and $b'_i$  
replaced by $a_i$ and $b_i$, respectively, and $R$ consists of the words of $\wh R$ with $\wh a_i$ and $\wh b_i$ replaced by $a_i$ and $b_i$, respectively. 
This implies the existence of the  required short surjectivity diagram. 
\end{proof} 

\begin{cory}\label{c:prod}
Let $F\map{i} X \to S^2$ be a symplectic Lefschetz fibration. Let $Y=\Si_e\times 
F$. Then $\pi_1 (X\#_FY)\simeq \pi_1(X)\times \pi_e$.
\end{cory}

\begin{proof}
Consider the diagram of \lemref{L:pi1}. Clearly, its top exact sequence splits, and therefore the bottom exact sequence splits by \propref{p:ssd}.
\end{proof}

\begin{cory}\label{c:split}
Assume that a group $\Ga$ admits a finite representation $p_{\Ga}: \pi_g\to \Ga 
$ such that 
$$
\pg^*: H^2(\Ga;\R)\to H^2(\pi_g;\R)
$$ is a non-zero homomorphism. Then the group $\Ga \times \pi_e$ is 
symplectically aspherical for all $e>0$.
\end{cory}

\begin{proof} Consider the symplectic Lefschetz fibration $F\to X \to S^2$ as in 
\propref{P:fac}. Let $Y=\Si_e\times F$. Then the manifold $X\#_FY$ is symplectically aspherical by \propref{p:sum}, and the claim follows from \coryref{c:prod}
\end{proof}


\section{Symplectically aspherical abelian groups}\label{S:sab}


In this section we give a complete description of symplectically aspherical 
abelian groups.

\begin{lemma}\label{l:epi} 
Let $\Ga $ be a finitely presentable abelian group such that $H^2(\Ga; \R)\neq 
0$. 
Then there exists a finite presentation $p_{\Ga}:\pi_g\to \Ga$ such that the map 
$$
p^*_{\Ga}: H^2(\Ga;\R) \to H^2(\pi_g;\R)
$$
is non-zero.
\end{lemma}

\begin{proof} First, it is easy to construct a finite presentation $\pi_h\to 
\Ga$ (and here there is no necessity to assume $\Ga$ abelian). Indeed, there are 
two finite presentations $F_r\to \Ga$ and $\pi_{2r} \to F_r$ where $F_r$ is the free 
group of $r$ generators. Now, the composite $f:\pi_{h}\to F_r\to \Ga, h=2r$ is a 
finite presentation.
Since $H^2(\Ga; \R)\neq 0$, we conclude that there exists a monomorphism $g: \Z^2\to \Ga$ such that the induced homomorphism  
$$
g^*:H^2(\Ga;\R)\to H^2(\Z^2;\R)
$$
 is non-trivial. Consider the map 
$$
\gf: \pi_{h}*\Z^2 \to \Ga,\quad \gf(a_1b_1\cdots a_nb_n)=f(a_1)g(b_1)\cdots 
f(a_n)g(a_n).
$$
We claim that it is a finite presentation of $\Ga$. Notice the following:
\begin{enumerate} 
\item If $f: A \to B$ and $g: B \to C$ are finite presentations then so is $gf: 
A \to C$.
\item If $A$ is a finitely presentable group then the  abelianization $\ab: A \to 
A_{\ab}$ is a finite presentation.
\end{enumerate}

Clearly, $\pi_{h}*\Z^2$ is a finitely presentable group. Now, since $\Ga$ is an 
abelian group, the epimorphism $\gf$ can be decomposed as 
$$
\CD
\gf: \pi_{h}*\Z^2 @>\ab >> (\pi_{h}*\Z^2)_{\ab} @>>> \Ga
\endCD
$$
where both maps are finite presentations (the last one because 
$(\pi_{h}*\Z^2)_{\ab}$ is a finitely generated abelian group). So, $\gf$ is a 
finite presentation by (1) and (2).

Now, there is a canonical map $u: \Si_{h+2}\to \Si_h\vee T^2$ (namely, we regard 
$\Si_{h+2}$ as the connected sum  $\Si_h\# T^2$ and pinch the circle which we 
glued the surfaces along). Clearly, $u_*: \pi_1(\Si_{h+2}) \to \pi_1(\Si_h \vee T^2)$ 
is a finite presentation (its kernel 
is generated as a normal subgroup by the pinched circle). Thus, by (1), the 
homomorphism
$$
\CD
\pi_{h+2}=\pi_1(\Si_{h+2}) @>u_*>> \pi_1(\Si_h)*\pi_1(T^2) @= \pi_h*\Z^2 @>\gf 
>> \Ga
\endCD
$$
is the desired finite presentation of $\Ga$.
\end{proof}

\begin{proof}[Proof of Theorem \ref{T:abel}] 
Let $T$ be a finite abelian group and 
$\Ga=\B Z^m\oplus T$ with $m\ge 4$. 
Let $A= \B Z^{m-2}\oplus T$. 
Then $H^2(A;\R)\ne 0$, and so, by \lemref{l:epi} and 
\corref{c:split} applied to the case 
$e=1$, we conclude that the group 
$\Ga=A\oplus \Z^2$ is symplectically aspherical.

Furthermore, the group $\B Z^3\oplus T$ is not symplectically 
aspherical because it is three dimensional 
(see Example \ref{E:3dim}).
Finally, suppose that the group $\B Z^2\oplus T=\pi_1(M)$ for $M$ closed symplectically 
aspherical. Since $H^i(\B Z^2\oplus T)=0$ for $i>2$, we conclude that
$M$ must be 2-dimensional by \cite[Proposition 2.3]{IKRT}). So, $M$ is a closed surface 
(in fact, the torus $\T^2$), and thus $T=0$. 
\end{proof} 

\begin{cory}\label{cor:dim}
If $m>3$, then every abelian group $\Z^m \oplus T$ with $T$ finite can be 
realized as the fundamental group of a closed symplectically aspherical manifold $M^{2n}$ 
with $4\le 2n \le m$ and cannot be realized as the fundamental group of a closed 
symplectically aspherical manifold $M^{2n}$ with $2n>m$.
\end{cory}

\begin{proof} Take $n$ with $4\le 2n \le m$. We know that there exists 
4-dimensional manifold $N$ with $\pi_1(N)=\Z^{m-2n+4}\oplus T$. Now, let 
$M=N\times \T^{2n-4}$.

The last claim follows from \cite[Proposition 2.3]{IKRT}.
\end{proof}

\begin{prop}\label{prop:pi2}
Let $N$ be a closed $4$-dimensional symplectically aspherical manifold such 
that $\pi_1(N)=\Z^4\oplus \Z/2$. Then $\pi_2(N)\ne 0$.
\end{prop}

Notice that Theorem~\ref{T:abel} guarantees the existence of such manifold $N$.

\begin{proof} Indeed, in case $\pi_2(N)=0$ we have the Hopf exact sequence 
$$
\CD
\pi_3(N) @>>>  H_3(N) @>>>  H_3(\pi_1(N)) @>>> 0.
\endCD
$$
But $H_3(N)=H^1(N)=\Hom (\B Z^4\oplus \Z/2,\, \Z)=\B Z^4$ by the Poincar\'e duality 
and the universal coefficient theorem, while 
$$
H_3(\B Z^4\oplus \Z/2)\supset H_3(\B Z^4) \oplus H_3(\B Z/2)=\B
Z^4\oplus \B Z/2.
$$ 
Hence, there are no epimorphisms $H_3(N) \to  H_3(\pi_1(N))$, and thus
$\pi_2(N)\ne 0$.
\end{proof}

Question 8.3.2 in \cite{IKRT} asks whether there exists a closed symplectically 
aspherical manifold 
$M$ with $\pi_1(M)=\B Z^4$ and $\pi_2(M)\ne 0$. Now we can answer affirmatively.

\begin{cory}\label{cor:z4}
There exists a closed $4$-dimensional symplectically aspherical manifold with 
$\pi_1(M)=\B Z^4$ and $\pi_2(M)\ne 0$. 
\end{cory} 

\begin{proof} Let $N$ be a manifold considered in \propref{prop:pi2} and $\omega$ be a symplectically aspherical form on $N$.
Let $p:M \to N$ be a two-sheeted covering with
$\pi_1(M)=\Z^4$. Then $p^*\omega$ is a symplectically aspherical form on a closed manifold $M$,
while $\pi_2(M)\ne 0$ and $\pi_1(M)=\B Z^4$.
\end{proof}

Question 8.3.1 in \cite{IKRT} asks whether every symplectically aspherical 
group can be realized as the fundamental group of a  closed symplectically aspherical 
manifold $M$ with $\pi_2(M)=0$.

\begin{cory}\label{cor:pi2} 
If $N$ is a  closed symplectically aspherical manifold such 
that  $\pi_1(N)=\Z^4\oplus \Z/2$, then $\pi_2(N)\ne 0$.
\end{cory}

\begin{proof} Indeed, because of \corref{cor:dim}, we must have $\dim N=4$. But 
then $\pi_2(N)\ne 0$ by \propref{prop:pi2}.
\end{proof}

\end{document}